\documentclass[reqno,11pt]{amsart}

\usepackage{latexsym,amssymb,amsthm,amsmath,amscd}

\usepackage{amssymb}
\usepackage{eucal}
\usepackage{amsmath}

\theoremstyle{plain}
\newtheorem{theorem}{Theorem}[section]
\newtheorem*{Theorem B}{Theorem B}
\newtheorem*{Theorem A}{Theorem A}

\newtheorem{definition}{Definition}[section]

\numberwithin{equation}{section}

\theoremstyle{remark}

 \numberwithin{equation}{section}

\begin{document}

\noindent  {\footnotesize J. Geom. Symmetry Phys. 43 (2017)
\vskip-.05in

\noindent {\it Accepted paper}}

\vskip.3in

\title[Covering maps and ideal embeddings]{Covering maps and ideal embeddings of compact homogeneous spaces}

\author[B.-Y. Chen]{Bang-Yen Chen}
\address{Department of Mathematics, Michigan State University, East Lansing, Michigan 48824--1027, U.S.A.}
\email{bychen@math.msu.edu}

\begin{abstract}
The notion of ideal embeddings was introduced in [B.-Y. Chen, {Strings of Riemannian invariants,  inequalities, ideal immersions  and their applications.} The Third Pacific Rim Geometry Conference (Seoul, 1996), 7--60,  Int. Press, Cambridge, MA, 1998]. Roughly speaking, an ideal embedding (or a best of living) is an isometrical embedding which receives the least possible amount of tension from the surrounding space at each point.

In this article, we study ideal embeddings of irreducible compact homogenous spaces in Euclidean spaces via covering maps. Our main result states that {if $\pi: M\to N$ is a covering map between two irreducible compact  homogeneous spaces and if $\lambda_1(M)\ne \lambda_1(N)$, then $N$  doesn't admit an ideal embedding in any Euclidean space; although $M$ could}.  \\[0.2cm]
 \textsl{MSC}: {53C30; 53C40; 53C42} \\
 \textsl{Keywords}: Ideal embedding; best way of living; compact homogeneous spaces;  first positive eigenvalue of Laplacian; covering map.
\end{abstract}

\subjclass[2000]{53C30; 53C35; 53C42}
\keywords{Ideal embedding; best way of living; compact homogeneous spaces;  first positive eigenvalue of Laplacian; covering map.}

\maketitle

\section{Introduction}

According to John F. Nash's embedding theorem \cite{nash}, every Riemannian manifold can be isometrically embedded in a Euclidean space with sufficiently large codimension. In other words, every Riemannian manifold can live in the Euclidean world if the codimension was sufficiently large. 

Related to Nash's theorem, my main question raised in \cite{c97,c98} is the following.

\vskip.1in
\noindent{\bf Main Question.} {\it Can a given Riemannian manifold live in a Euclidean world ideally? 

More precisely,  can a given Riemannian manifold be isometrically embedded in a Euclidean space in such way that it receives the least possible tension from the surrounding space at each point?}
\vskip.1in

It is well-known that the mean curvature vector field of a submanifold is exactly the tension field for an isometric immersion of a Riemannian manifold in another Riemannian manifold.

A major concern of my main question is whether we can determine when a Riemannian manifold can be embedded in a Euclidean space ideally.
Several answers to this question had been obtained (cf. \cite{c97,c98,c00a,c08,book,C13b} among others). 

For a compact Riemannian manifold $M$, we denote the first positive eigenvalue of the Laplacian $\Delta_M$ of $M$ by $\lambda_1(M)$. By a {\it covering map} $\pi: M\to N$ between two Riemannian manifolds we mean {\it a covering map which is isometric}.

The following result on $\lambda_1$ was proved by K. Yoshiji in \cite{Yo00}.
 
 \begin{theorem} Every non-orientable compact manifold $M$, except the real projective plane ${\mathbb R}P^2$,  admits a Riemannian metric $g$ for which the first eigenvalue coincides with that of its Riemannian double cover $\hat M$, i.e., one has $\lambda_1(M)=\lambda_1(\hat M)$ with respect to the metric $g$ and its covering metric $\hat g$ on $\hat M$.
  \end{theorem}

A Riemannian manifold is said to be {\it homogeneous} if the group of isometries of $M$ acts transitively on $M$ (cf. \cite{He78,KN63}). It is well-known that every homogeneous Riemannian manifold is complete (cf. \cite[Theorem 4.6]{KN63}). 
A compact  homogeneous Riemannian manifold with irreducible isotropy action is simply called an {\it irreducible compact  homogeneous space}.

Homogeneous spaces are important. From the point of view of the {\it Erlangen program}, in a homogeneous space one may understand that ``all points are the same''. 
For example, Euclidean space and projective space are in natural ways homogeneous spaces for their respective symmetry groups. The same is true of the models found of non-Euclidean geometry of constant curvature, such as hyperbolic space.

In this article, we will investigate the Main Question further. In 
particular, by studying covering maps between irreducible 
homogeneous spaces, we are able to obtain a solution to the 
Main Question for irreducible compact homogeneous spaces.

The main result of this article is the following. 
\vskip.1in

\noindent {\bf Main Theorem.} {\it Let $\pi: M\to N$ be a covering map between two irreducible compact  homogeneous spaces. If $\lambda_1(M)\ne \lambda_1(N)$, then $N$  doesn't admit an ideal embedding in any Euclidean space, regardless of codimension}.  

\vskip.1in
 
In the last section, we give a simple example to illustrate that, under the hypothesis of the Main Theorem, $M$ may admit an ideal embedding in some Euclidean space; although $N$ cannot.

\section{$\delta$-invariants and fundamental inequality}

In order to define ideal embedding, we need to recall the notion of $\delta$-invariants (also known as the {\it Chen invariants}\/) and the fundamental inequality of Euclidean submanifolds.

 Let $M$ be a  Riemannian $n$-manifold. Let $K(\pi)$ denote the sectional curvature of $M$ associated with a plane section $\pi\subset T_pM$, $p\in M$. For a given orthonormal basis $e_1,\ldots,e_n$ of the tangent space $T_pM$, the scalar curvature $\tau$ at $p$ is defined to be $$\tau(p)=\sum_{i<j} K(e_i\wedge e_j).\notag$$

Let $L$ be a subspace of $T_pM$  of dimension $r\geq 2$  and let $\{e_1,\ldots,e_r\}$ be an orthonormal basis of $L$. We define the scalar curvature $\tau(L)$ of $L$ by  $$\tau(L)=\sum_{\alpha<\beta} K(e_\alpha\wedge e_\beta),\quad 1\leq \alpha,\beta\leq r.$$

Given an integer $k\geq 1$, we denote by ${\mathcal S}(n,k)$ the finite set  consisting of unordered $k$-tuples $(n_1,\ldots,n_k)$ of integers $\geq 2$ satisfying  $n_1< n$ and $n_1+\cdots+n_k\leq n$. We put ${\mathcal S}(n)=\cup_{k\geq 1}{\mathcal S}(n,k)$. 

 For each $k$-tuple $(n_1,\ldots,n_k)\in {\mathcal S}(n)$, Chen defined the {\it $\delta$-invariant} $\delta{(n_1,\ldots,n_k)}$ in \cite{c98,c2000,book} as 
\begin{align} \delta(n_1,\ldots,n_k)(p)=\tau(p)-\inf\{\tau(L_1)+\cdots+\tau(L_k)\}, \end{align}
where $L_1,\ldots,L_k$ run over all $k$ mutually orthogonal subspaces of $T_pM$ such that  $\dim L_j=n_j,\, j=1,\ldots,k$. 
In particular, we have
\vskip.05in 

$\delta(\emptyset)=\tau$ \ \   ($k=0$, the trivial $\delta$-invariant),
\vskip.05in

$\delta{(2)}=\tau-\inf K,$ where $K$ is the sectional curvature,

\vskip.05in
$\delta(n-1)(p)=\max Ric(p)$,

where $Ric$ is the Ricci curvature of $M$.
\vskip.1in
We shall point out that the invariant $\delta(2)$ was introduced earlier in \cite{c93,c95}.

The $\delta$-curvatures are very different in nature from the ``classical'' scalar and Ricci curvatures; simply due to the fact that both scalar and Ricci curvatures are the ``total sum'' of sectional curvatures on a Riemannian manifold. In contrast, the $\delta$-curvature invariants  are obtained  from the scalar curvature by throwing away a certain amount of sectional curvatures. 
For history and motivation on $\delta$-invariants, see author's surveys \cite{c08,book,C13b}.

For an isometric immersion of a Riemannian $n$-manifold $M$ into the Euclidean $m$-space $\mathbb E^m$. Let $h$ and $\overrightarrow{H}$ denote the second fundamental form and the mean curvature vector of $M$ in $\mathbb E^m$, respectively. Then $h$ and $\overrightarrow{H}$ are defined respectively by
\begin{align}\label{2.2} &h(X,Y)=\tilde \nabla_XY-\nabla_X Y,
\\&\label{2.3}  \overrightarrow{H}=\frac{1}{n} \,{\rm trace},\end{align} for vector fields $X$ and $Y$ tangent to $M$,
where $\tilde \nabla$ and $\nabla$ denote the Levi-Civita connection on $\mathbb E^m$ and on $M$, respectively. We put $H=|\overrightarrow{H}|$.

The fundamental inequality for  Euclidean submanifolds obtained in \cite{c98,c00a,book} is the following sharp inequality.

\begin{theorem} \label{T:2.1} For any isometric immersion of a Riemannian $n$-manifold $M$ into a Euclidean $m$-space $\mathbb E^m$ and for any $k$-tuple $(n_1,\ldots,n_k)\in{\mathcal S}(n)$, we have
\begin{align}\label{2.4}\delta{(n_1,\ldots,n_k)} \leq  c(n_1,\ldots,n_k)H^2,
\end{align}
where $c(n_1,\ldots,n_k)$ is the positive real number defined by
\begin{align}\label{2.5}c(n_1,\ldots,n_k)= {{n^2(n+k-1-\sum n_j)}\over{2(n+k-\sum n_j)}}.
\end{align}
\end{theorem}

\section{Ideal embeddings - best ways of living}

The fundamental inequality \eqref{2.4} in Theorem \ref{T:2.1} provides us the prime controls on the most important extrinsic curvature, the squared mean curvature $H^2$,  by the $\delta$-invariant $\delta(n_1,\ldots,n_k)$  of the Riemannian manifold $M$.

If we put 
\begin{align}&
\label{3.1}\hat\Delta_0(M)=\max\,\{\Delta(n_1,\ldots,n_k):(n_1,\ldots,n_k)\in {\mathcal S}(n)\} \end{align} 
with \begin{align}\Delta(n_1,\ldots,n_k)=\frac{\delta(n_1,\ldots,n_k)}{c(n_1,\ldots,n_k)},\end{align}
then, for any isometric immersion  $x:M\to  \mathbb E^m$, Theorem \ref{T:2.1} yields
\begin{align}\label{3.2} H^2\geq \hat\Delta_0(M)\end{align} 
at each point $p\in M$.

Inequality \eqref{3.2} allowed us to introduce the following notion of ideal immersions in \cite{c98,c00a,book}.

\begin{definition} {\rm An isometric immersion  of a Riemannian $n$-manifold $M$ in $\mathbb E^m$ is called an {\it ideal immersion\/} if it satisfies the equality case of \eqref{3.2} identically.} \end{definition}

A {\it Maximum Principle} on ideal immersions from \cite{c98} states that
if an isometric immersion $\phi:N\to \mathbb E^m$ of a Riemannian $n$-manifold $M$ into $\mathbb E^m$ satisfies the equality sign of \eqref{3.1} for a given $k$-tuple $(n_1,\ldots,n_k)\in {\mathcal S}(n)$, then it is an ideal immersion automatically.

\vskip.1in
\noindent {\bf Physical Interpretation of Ideal Immersions}.  An isometric immersion  $x:M\to \mathbb E^m$ is an {\it ideal immersion} means that the submanifold $M$  receives the least possible amount of tension  (given by $\hat\Delta_0(M)$) from the surrounding space at each point on $M$. This is due to \eqref{3.2} and the well-known fact that the mean curvature vector field is exactly the tension field for an isometric immersion of a Riemannian manifold in another Riemannian manifold; thus  the squared mean curvature at each point on the submanifold simply measures the amount of tension the submanifold is receiving from the surrounding space at that point. 

For this reason,  an ideal embedding of a Riemannian manifold $M$ in a Euclidean space is also known as a {\it best way of living} of $M$ (see  \cite{c97,c98,book}).
\vskip.1in

It was shown in \cite{C13a} that some Riemannian manifolds may admit more than one ideal embedding in a Euclidean space.

A major problem on ideal embeddings is to determine whether a given Riemannian manifold  admits or doesn't admit an ideal embedding in a Euclidean space.

\section{Proof of the Main Theorem}

 Let $\pi: M\to N$ be a covering map from a compact Riemannian manifold onto another. Suppose that $f\in C^{\infty}(N)$ is an eigenvalue function of the Laplacian $\Delta_{N}$ associated with the first positive eigenvalue $\lambda_1(N)$ of $N$. Then we have
\begin{align}\label{4.1}  \Delta_M (f\circ \pi)=\lambda_1(N) f \circ\pi.\end{align}
It follows from \eqref{4.1} that $\lambda_1(N)$ is an eigenvalue of the Laplacian $\Delta_M$ of $M$. Thus we  have
\begin{align}\label{4.2} \lambda_1(N)\geq \lambda_1(M).\end{align}

Now, suppose that $\pi: M\to N$ is a covering map from a compact irreducible homogeneous space $M$ onto another  irreducible compact homogeneous space $N$ and suppose that  $\lambda_1(M)\ne \lambda_1(N)$. Then, by combining $\lambda_1(M)\ne \lambda_1(N)$ with \eqref{4.2} we obtain
\begin{align}\label{4.3}\lambda_1(N)> \lambda_1(M).\end{align}

On the other hand,  from Theorem 6.3 of \cite{c00a} or Theorem 14.6 of \cite{book} we know that the first positive eigenvalue $\lambda_1(M)$ of the irreducible homogeneous space $M$ satisfies
\begin{align}\label{4.4}\lambda_1(M)\geq n \hat \Delta_0(M).\end{align}
Since $M$ is assumed to be a compact homogeneous Riemannian manifold,  $\hat\Delta_0(M)$ defined in \eqref{3.1} is a constant  on $M$. Because the covering map
 $\pi: M\to  N$ preserves the metric tensors $g_M$ and $g_{N}$ of $M$ and $N$, i.e., $\pi^*g_{N}=g_M$, we know that the invariant $\hat\Delta_0(N)$ of $N$ must equals to  $\hat\Delta_0(M)$ on $M$. Therefore, after combining \eqref{4.3} and \eqref{4.4} we obtain 
  \begin{align}\label{4.5}\lambda_1(N)>n \Delta_0(N).\end{align}
 
 From Theorem 6.6 of \cite{c00a} or Theorem 14.7 of my book \cite{book} we also know that an irreducible compact  homogeneous space $M$ admits an ideal immersion into some Euclidean space if and only if it satisfies 
 \begin{align}\label{4.6}\lambda_1(M)= n \hat \Delta_0(M).\end{align}
 Consequently,  inequality \eqref{4.5} together with \cite[Theorem 14.7]{book} imply that $N$  does not admit an ideal embedding into Euclidean spaces, regardless of codimension. This completes the proof of the Main Theorem.

\section{An example}

Now, we provide a simple example to illustrate that under the hypothesis of the Main Theorem, $M$ may admit an ideal embedding in some Euclidean space; although $N$ cannot.

 Let $S^n(1)$ and ${\mathbb R}P^n(1)$ denote the $n$-sphere and the real projective $n$-space of constant sectional curvature one. It is well-known that there is a two-fold covering map $\pi:S^n(1)\to {\mathbb R}P^n(1)$ which carries each pair of antipodal points on $S^n(1)$ to a single point in ${\mathbb R}P^n(1)$. 

Since both $$S^n(1)=SO(n+1)/SO(n)\;\; {\rm and} \;\; {\mathbb R}P^n(1)=SO(n+1)/SO(n)\times \{\pm 1\}$$ are irreducible compact homogeneous spaces and moreover we have (cf. e.g. \cite{BGM71} and \cite[page 67]{c15})
 $$n=\lambda_1(S^n(1)) \ne \lambda_1({\mathbb R}P^n(1))=2(n+1),$$  the Main Theorem thus implies that $RP^n(1)$ never admits a best way of living in any Euclidean space regardless of codimension.  

On the other hand, it is easy to verify that the inclusion map of $S^n(1)\subset \mathbb E^{n+1}$ is an ideal embedding of $S^n(1)$ in $\mathbb E^{n+1}$.

 Bang-Yen Chen \\
Department of Mathematics \\
Michigan State University\\
East Lansing, Michigan 48824-1027\\ USA\\
{\it E-mail address}: {\tt bychen@math.msu.edu}

\label{last}
\end{document}